\documentclass[11pt]{article}
\usepackage{amsmath}
\usepackage{latexsym}
\usepackage{amssymb}
\usepackage{theorem}
\usepackage{amsfonts}

\usepackage{graphicx}
\usepackage{amsmath}
\usepackage{amssymb}
\usepackage{bm}
\usepackage[all]{xypic}
\usepackage[all]{xy}

\theoremstyle{plain} \theoremstyle{} {\theorembodyfont{\itshape}
\newtheorem{thm}{Theorem}[section]}
\theoremstyle{plain} \theoremstyle{} {\theorembodyfont{\itshape}
\newtheorem{lem}[thm]{Lemma}}
\theoremstyle{plain} \theoremstyle{} {\theorembodyfont{\itshape}
\newtheorem{prop}[thm]{Proposition}}
\theoremstyle{plain} \theoremstyle{} {\theorembodyfont{\itshape}
}
\theoremstyle{plain} \theoremstyle{} {\theorembodyfont{\itshape}
\newtheorem{cor}[thm]{Corollary}}
\theoremstyle{plain} \theoremstyle{} {\theorembodyfont{\itshape}
}
\theoremstyle{plain} \theoremstyle{} {\theorembodyfont{\itshape}
}
\theoremstyle{plain} \theoremstyle{} {\theorembodyfont{\rmfamily}
}
\theoremstyle{plain} \theoremstyle{} {\theorembodyfont{\rmfamily}
}
\theoremstyle{plain} \theoremstyle{} {\theorembodyfont{\rmfamily}
}
\theoremstyle{plain} \theoremstyle{} {\theorembodyfont{\rmfamily}
}
\theoremstyle{plain} \theoremstyle{} {\theorembodyfont{\rmfamily}
\newtheorem{rks}[thm]{Remarks}}
\theoremstyle{plain} \theoremstyle{} {\theorembodyfont{\rmfamily}
}
\theoremstyle{plain} \theoremstyle{} {\theorembodyfont{\rmfamily}
}
\theoremstyle{plain} \theoremstyle{} {\theorembodyfont{\rmfamily}
}
\theoremstyle{plain} \theoremstyle{} {\theorembodyfont{\rmfamily}
\newtheorem{df}[thm]{Definition}}

\theoremstyle{plain} \theoremstyle{} {\theorembodyfont{\itshape}
}
\theoremstyle{plain} \theoremstyle{} {\theorembodyfont{\itshape}
}
\theoremstyle{plain} \theoremstyle{} {\theorembodyfont{\itshape}
}
\theoremstyle{plain} \theoremstyle{} {\theorembodyfont{\rmfamily}
}

\begin{document}
\input xy
\xyoption{all}

\newcommand{\Hom}{\mathrm{Hom}}

\newcommand{\End}{\mathrm{End}}
\newcommand{\Ext}{\mathrm{Ext}}
\newcommand{\rad}{\mathrm{rad}}
\newcommand{\Aut}{\mathrm{Aut}}
\newcommand{\ind}{\mathrm{ind}}

\newcommand{\moa}{\mbox{mod}{A}}
\newcommand{\Z}{{\bf Z}}
\newcommand{\C}{{\bf C}}
\newcommand{\Q}{{\bf Q}}
\newcommand{\N}{{\bf N}}

\newcommand{\g}{\bf{g}}
\newcommand{\n}{\bf{n}}
\newcommand{\h}{\bf{h}}

\newcommand{\al}{\alpha}
\newcommand{\be}{\beta}
\newcommand{\ga}{\Gamma}
\newcommand{\si}{\sigma_i}
\newcommand{\sip}{\sigma_i^+}
\newcommand{\de}{\Delta}
\newcommand{\la}{\longleftarrow}
\newcommand{\ra}{\longrightarrow}
\newcommand{\xra}{\xrightarrow}
\newcommand{\va}{\varepsilon}
\newcommand{\ad}{\mbox{ad}}
\newcommand{\vlmta}{V_{\lambda_{m-t-a}}}
\newcommand{\mc}{\mathcal{C}}
\newcommand{\md}{\mathcal{D}}
\newcommand{\mh}{\mathcal{H}}
\newcommand{\mr}{\mathcal{R}}
\newcommand{\ma}{\mathcal{A}}
\newcommand{\ml}{\mathcal{L}}
\newcommand{\mt}{\mathcal{T}}
\newcommand{\me}{\mathcal{E}}
\newcommand{\ms}{\mathcal{S}}
\newcommand{\mk}{\mathcal{K}}
\newcommand{\mpi}{\mathcal{P}}
\newcommand{\zr}{\mathbb{Z}_2}
\newcommand{\tauni}{\tau^{-1}}
\newcommand{\taun}{\tau^{n}}
\newcommand{\taunin}{\tau^{-n}}
\newcommand{\homc}{\Hom_{\mathcal{C}}}
\newcommand{\homh}{\Hom_{\mathcal{H}}}
\newcommand{\raa}{\overrightarrow{a}}

\title{{\Large \textbf{ The Bridgeland's Ringel-Hall algebra associated to an algebra with global dimension at most two
}}\thanks{Supported partily by the National Natural Science Foundation of China (Grant No. 11126319
) and the National 973 Programs (Grant No. 2011CB808003).}}
\author{ {\normalsize \textbf{Shengfei Geng}}\\
{\small\sl Department of Mathematics, Sichuan University, Chengdu
610064, China} \\
{\small\sl  E-mails: genshengfei@{scu.edu.cn}
}\\ {\normalsize  \textbf{Liangang Peng}}\\
{\small\sl Department of Mathematics, Sichuan University, Chengdu
610064, China} \\
{\small\sl  E-mails: penglg@scu.edu.cn}}
\date{}
\maketitle

 \maketitle \noindent \textbf{Abstract.} For any finite dimensional associative algebra with global dimension $\leq 2$, we show that there is an embedding from the twisted Ringel-Hall algebra to the Brigeland's Ringel-Hall algebra. In particular, this result is true for tilted algebras and canonical algebras.
 \vspace{0.3cm}


\vspace{0.3cm}

\parskip 3pt

\section{\large  Introduction}

For any finite dimensional  semi-simple complex Lie algebra,  there is a famous theorem saying that its positive roots correspond bijectively to the isomorphism classes of all (finitely dimensional) indecomposable modules of the  corresponding hereditary algebra, which was proved by Gabriel  \cite{Gab} for ADE type and then extended by Dlab-Ringel \cite{DR} for any type. This result gives a realization of the positive root system of the semi-simple complex Lie algebra from the modules of the hereditary algebra.

To realize the multiplication of the semi-simple complex Lie algebra from the module category of the hereditary algebra, Ringel  \cite{Ri2} introduced a Hall algebra from any abelian category  with finite morphism spaces, which generalizes Hall's definition from $p$-groups \cite{Ha}. 
 He \cite{Ri3} showed that the  positive part of the semi-simple complex Lie algebra can be realized via the Hall algebra from the module category of the hereditary algebra. Moreover, he  \cite{Ri4}
proved that the  positive part of  the   quantized enveloping algebra of the semi-simple complex Lie algebra is isomorphic to the Hall algebra in some twisted form from the module category of the hereditary algebra. Later,  Green  \cite {Gr} extended this result to any symmetrizable Kac-Moody Lie algebra case by showing that there is a natural co-multiplication in the Hall algebra of any hereditary algebra.

So a natural question is how one can use the representation theory of hereditary algebras to realize whole but not only the positive part of any symmetrizable Kac-Moody Lie algebra and its quantized enveloping algebra. There have been many attempts to approach such a question. One way was to consider the Drinfeld double of the Hall algebra of any hereditary algebra constructed  by Xiao in \cite{X},  and in this way he realized the full quantized enveloping algebra. Second way was given by Peng-Xiao in \cite{PX1} and \cite{PX2} who used the similar method as Ringel's to construct the Ringel-Hall Lie algebra from any root category, which is the orbit category of all shift square orbits in the derived category of a hereditary algebra. In this way they realized the full symmetrizable Kac-Moody Lie algebra. Third way was introduced by To\"{e}n in \cite{T} who constructed the derived  Ringel-Hall algebra from the derived category. Recently, Bridgeland  \cite{Br} gave a more clever way to consider directly the $\mathbb{Z}_2$-graded complexes (or 2-cycle complexes as called in \cite{PX1},  or 2-periodic complexes as in some literature) of projective modules of an algebra with finite global dimension to construct the Ringel-Hall algebra. He showed that the full quantized enveloping algebra can be also realized by his Ringel-Hall algebra. Later Yanagida  \cite{Y} proved that in hereditary algebra case the Bridgeland's Ringel-Hall algebra is isomorphic to the Drinfeld double of the Hall algebra. Inspired by work of Bridgeland, Grosky recently in \cite{Gro} constructed the semi-derived Ringel-Hall algebras from  complex categories or $\mathbb{Z}_2$-graded complex categories. In each way above there also were some further researches, see    \cite{Ka}, \cite{SVdB1}, \cite{SVdB2}, \cite{PT}, \cite{Cramer},\cite{Hubery}, \cite{LP},\cite{Sch},\cite{XX}, and so on.

In this paper we only consider Bridgeland's Ringel-Hall algebras. For any hereditary algebra, Bridgeland \cite{Br} has proven that the twisted Ringel-Hall algebra can be naturally embedded in the Bridgeland's Ringel-Hall algebra. We extend this result to the case of any algebra  with global dimension $\leq 2$. Some special interesting consequences should be for tilted algebras and canonical algebras. Such two kinds of algebras are very important  and they both have global dimension $\leq 2$.

\noindent{ \bf Notation}. We fix a field $k = \mathbb{F}_q$ with $q$ elements, and set $t=\sqrt{q}$. We write $[A]$ for the
isomorphism class of  $A\in\mathcal{A}$ and $\mathrm{Iso}(\mathcal{A})$ for the set of isomorphism classes of a small category $\mathcal{A}$. The symbol $|S|$ denotes the number of elements of a finite set $S$.

\section{Preliminaries }

In this section, we give basic definitions and properties of  Ringel-Hall algebra, following   \cite{Ri4}, \cite{Br}, and \cite{Sch}.

In this section,  $\ma$ denotes an abelian category and satisfy the following conditions:

$(a)$ $\ma$ is essentially small, with finite morphism spaces,

$(b)$ $\ma$ is linear over $k = F_q$,

$(c)$ $\ma$ is of finite global dimension and has enough projectives.

\

 Given objects $A,B,C \in\ma$, define $\Ext^1_{\ma}(A,C)_B \subset\Ext^1_{\ma} (A,C)$ to be the subset parameterising extensions with middle term isomorphic to $B$.
 \begin{df}\label{def of hall alg}
 The Ringel-Hall algebra $\mh(\ma)$ is  the vector space over $\mathbb{C}$ with basis indexed by elements $[A]\in \mathrm{Iso}(\ma)$, and with associative multiplication defined by$$[A]\diamond[C] = \sum_{[B]\in \mathrm{Iso}(\ma)}\frac{ |\Ext^1_{\ma}(A,C)_B|}{ |\Hom_{\ma}(A,C)|}[B].$$ The unit is [0].
\end{df}

 Let $K(\mathcal{A})$ denote the Grothendieck group of $\mathcal{A}$. We write $\hat{A} \in K(\mathcal{A})$ for the class of an object $A \in \mathcal{A}$. Let $K_{\geq 0}(\mathcal{A}) \subset K(\mathcal{}A)$ be the subset consisting of these classes. For objects $A,B\in\mathcal{ A}$, define $$\langle  \hat{A},\hat{B}\rangle= \sum_{i\in \mathbb{Z}}(-1)^i\dim_k \Ext^i_\mathcal{A}(A,B).$$\noindent The sum is finite by our assumptions on $\mathcal{A}$, and descends to give a bilinear form $$\langle -,-\rangle: K(\mathcal{A})\times K(\mathcal{A})\rightarrow \mathbb{Z}$$
\noindent known as the Euler form. We also consider the symmetrized form $$(-,-): K(\mathcal{A})\times K(\mathcal{A})\ra \mathbb{Z},$$ \noindent defined by $(\alpha,\beta)=\langle \alpha,\beta\rangle+\langle \beta,\alpha\rangle$.

 \begin{df}
 The twisted Ringel-Hall algebra $\mathcal{H}_{tw}(\mathcal{A})$ is the same vector space as $\mh(\ma)$ but with twisted multiplication defined by $$[A]\ast [C]=t^{\langle \hat{A},\hat{C}\rangle}\cdot [A]\diamond [C].$$

\end{df}

 Let $\mc_{2}(\ma)$ be the abelian category of $\zr$-graded complexes in $\ma$. 
  The objects of this category consist of diagrams
$$ \xymatrix{
M_1\ar[r]<1mm>^{d_1}&M_0\ar[l]^{d_0}},\ \ \ \ d_{i+1}\circ d_i=0
$$
\noindent A morphism $s_{\bullet}:$ $M_{\bullet}\ra N_{\bullet}$
consists of a diagram
$$\xymatrix{&M_1\ar[d]_{s_1}\ar[r]<1mm>^{d_1}&M_0\ar[d]^{s_0}\ar[l]^{d_0}
\\&N_1\ar[r]<1mm>^{d'_1}&N_0\ar[l]^{d'_0}}$$
\noindent with $s_{i+1}\circ d_i=d'_i\circ s_i$. Two morphisms $s_{\bullet},t_{\bullet}: M_{\bullet}\ra N_{\bullet}$
are said to be homotopic if there are morphisms
$h_i : M_i \ra N_{i+1}$ such that
$$t_i-s_i = d'_{i+1}\circ h_i + h_{i+1}\circ d_i.$$
For an object $M_{\bullet}\in \mc_{2}(\ma)$,  we define its class in the $K$-group by
$$\hat{M_{\bullet}}:=\hat{M_0}-\hat{M_1}\in K(\ma).$$

Denote by $\mk_{2}(\ma)$ the homotopy category obtained from $\mc_{2}(\ma)$ by identifying homotopic morphisms. Let us also denote by $\mc_{2}(\mpi)\subset \mc_{2}(\ma)$
 the full subcategory whose objects are complexes of projectives  in $\ma$.
The shift functor [1] of complexes induces an involution $$\xymatrix{\mc_{2}(\ma)\ar[r]^{*}&\mc_{2}(\ma)\ar[l]}.$$
This involution shifts the grading and changes the sign of the differential as follows:
$$M_{\bullet}=(\xymatrix{M_1\ar[r]<1mm>^{d_1}&M_0\ar[l]^{d_0}})\xymatrix{\ar[r]^{*}&\ar[l]}M^{*}_{\bullet}=(\xymatrix{M_0\ar[r]<1mm>^{-d_0}&M_1\ar[l]^{-d_1}})$$

 Let $\md^b(\ma)$ denote the  bounded derived category of $\ma$, with its shift functor [1]. Let $R_2(A) = \md^b(\ma)/[2]$ be the orbit category, also known as the root category of $\ma$.
Since $\ma$ is assumed to be of finite global dimension with enough projectives, the category $\md^b(\ma)$ is equivalent to the ($\mathbb{Z}$-graded) bounded homotopy category $\mk^b(\mpi)$. Thus we can equally well define $R_2(\ma)$ as the orbit category of $\mk^b(\mpi)$.

\begin{lem} $\mathrm{(\cite{PX1},(\cite{Br})}$\label{rootcat and homolic cat} There is a fully faithful functor $D: R_2(\mathcal{A})\ra \mk_2(\mpi)$ sending a $\mathbb{Z}$-graded complex of projectives $(P_i)_{i\in \mathbb{Z}}$ to the $\mathbb{Z}_2$-graded complex $$\bigoplus_{i\in\mathbb{ Z}} P_{2i+1}\leftrightarrows \bigoplus_{i\in\mathbb{ Z}} P_{2i}.$$ \hfill{$\Box$}
\end{lem}
The functor $D$ is an equivalent when $\ma$ is hereditary (see \cite{PX1}).
\begin{lem} $\mathrm{(\cite{Br})}$\label{ext and hom over 2-} For $M_{\bullet},N_{\bullet}\in \mc_2(\mpi) $, we have
$$\Ext^1_{\mc_2({\ma})}(N_{\bullet},M_{\bullet})\cong \Hom_{\mk_2(\ma)}(N_{\bullet},M_{\bullet}^{*}).$$ \hfill{$\Box$}\end{lem}
A complex $M_{\bullet}\in \mc_2(\ma)$ is called acyclic if $H_{*}(M_{\bullet})=0.$ To each object $P \in \mpi$, we can attach acyclic complexes
$$K_{P_{}}=(\xymatrix{P\ar[r]<1mm>^{1}&P\ar[l]^0}), \ \ \ K^*_{P_{}}=(\xymatrix{P\ar[r]<1mm>^{0}&P\ar[l]^{1}}).$$

\begin{lem} $\mathrm{(\cite{Br})}$\label{accyci decomposition} For each acyclic complex of projectives $M_{\bullet}\in \mc_2(\mpi)$, there are objects $P,Q\in \mpi$, unique up to isomorphism, such that $M_{\bullet}\cong K_{P_{}}\oplus K_{Q_{}}^*.$\hfill{$\Box$}
\end{lem}
\

Let $\mh(\mc_2(\ma))$ be the Ringel-Hall algebra of the abelian category $\mc_2(\ma)$ defined in \ref{def of hall alg},
$\mh(\mc_2(\mpi))\subset \mh(\mc_2(\ma))$ be the subspace spanned by complexes of projective objects. Define $\mh_{tw}(\mc_2(\mpi))$ to be the same vector space as $\mh(\mc_2(\mpi))$ with the twisted multiplication $$[M_{\bullet}]\ast[ N_{\bullet}]:=t^{\langle M_0,N_0\rangle+\langle M_1,N_1\rangle}\cdot [M_{\bullet}]\diamond[N_{\bullet}],$$ then $\mh_{tw}(\mc_2(\mpi))$ is an associative algebra.

We have the following simple relations satisfied by the acyclic complexes $K_{P_{}},$
\begin{lem}$\mathrm{(\cite{Br})}$ \label{fact gongshi} For any object $P\in \mpi$ and any complex $M_{\bullet}\in \mc_2(\mpi)$, we have the following relations in $\mh_{tw}(\mc_2(\mpi)):$
     \begin{equation}[K_P]\ast [M_{\bullet}]=t^{\langle \hat{P},\hat{M_{\bullet}}\rangle}\cdot [K_P\oplus M_{\bullet}],\ \  [M_{\bullet}]\ast [K_P]=t^{-\langle \hat{M_{\bullet}},\hat{P}\rangle}\cdot[K_P\oplus M_{\bullet}]\label{fact gongshi1}\end{equation}
     \begin{equation} \label{fact gongshi jiaohuan}[K_P]\ast [M_{\bullet}]=t^{(\hat{P},\hat{M_{\bullet}})}\cdot[ M_{\bullet}]\ast[K_P],\ \  [K_P^*]\ast [M_{\bullet}]=t^{-(\hat{P},\hat{M_{\bullet}})}\cdot[ M_{\bullet}]\ast[K_P^*]\end{equation}
In particular, for $P,Q\in \mpi$, we have
     \begin{equation}[K_P]\ast[K_Q]=[K_P\oplus K_Q],\ \ [K_P]\ast[K_Q^*]=[K_P\oplus K_Q^*].\end{equation}\hfill{$\Box$}
\end{lem}

 \begin{lem} \label{fact gongshi2}
     \begin{equation}[K_P^*]\ast [M_{\bullet}]=t^{-\langle \hat{P},\hat{M_{\bullet}}\rangle}\cdot [K_P^*\oplus M_{\bullet}],\ \  [M_{\bullet}]\ast [K_P^*]=t^{\langle \hat{M_{\bullet}},\hat{P}\rangle}\cdot[K_P^*\oplus M_{\bullet}].\end{equation}
\end{lem}
{\bf Proof:} The proof is similar to the proof of the equations (\ref{fact gongshi1})(see \cite{Br}). \hfill{$\Box$}

\begin{df} The localized Ringel-Hall algebra $\md\mh(\ma)$ is the localization of $\mh_{tw}(\mc_2(\mpi))$ with respect to the elements $[M_{\bullet}]$ corresponding to acyclic complexes $M_{\bullet}:$
 $$\md\mh(\ma):=\mh_{tw}(\mc_2(\mpi))[[M_{\bullet}]^{-1}|H_{*}(M_{\bullet})=0]$$\end{df}
As explained in \cite{Br}, this is the same as localizing by the elements $[K_P]$ and $[K_P^*]$ for all objects $P\in \mpi$.

Write $\alpha\in K(\ma)$ in the form $\alpha=\hat{P}-\hat{Q}$ for objects $P,Q \in \mpi$, and set $K_{\alpha}=[K_{{P}}]\ast [K_{{Q}}]^{-1}, K_{\alpha}^*=[K_{{P}}^*]\ast [K^{*}_{{Q}}]^{-1}$. Still as explained in \cite{Br},   the equations  (\ref{fact gongshi jiaohuan})  continue to hold with the elements $[K_P]$ and $[K_P^*]$ replaced by $K_{\alpha}$ and $K^*_{\alpha}$ for arbitrary $\alpha \in K(\ma)$.

\begin{df} The reduced localized Ringel-Hall algebra is given by setting $[M_{\bullet}] = 1$ whenever $M_{\bullet}$ is an acyclic complex, invariant under the shift functor. In symbols
$$\mathcal{D}\mathcal{H}_{red}(\ma) =\mathcal{D}\mathcal{H}(\ma)/([M_{\bullet}]-1 : H_{*}(M_{\bullet}) = 0, M_{\bullet}\cong M_{\bullet}^*).$$
By Lemma \ref{accyci decomposition}, this is the same as setting $[K_P]\ast[K_P^*]=1$ for all $P\in \mathcal{P}.$
\end{df}

\section{Main results}
 In this section, let $\mathcal{A}=\mathrm{mod} B$ be the category of finitely generated right $B$-modules,
 where $B$ is a finite-dimensional  algebra with global dimension $\leq  2$. Hence, every object $A \in \mathcal{A}$
has a
projective resolution of the form
\begin{eqnarray}\label{minimal resolution} {0\ra P_2 \xra{a_2} P_1\xra{a_1} P_0\xra {a_0}A\ra 0.}\end{eqnarray}

The following is well known.

\begin{lem}
Any resolution (\ref{minimal resolution}) is isomorphic to a resolution of the form
\begin{eqnarray}\label{not minimal resolution}  0 \ra P'_2\oplus R_1 \xra{{\left(\begin{array}{cc}a'_2&0 \\0&0\\0&1\end{array}\right)}} P'_1\oplus R_0\oplus R_1\xra{{\left(\begin{array}{ccc}a'_1&0&0 \\0&1&0\end{array}\right)}} P'_0\oplus R_0\xra {(a'_0,0)}A\ra 0\end{eqnarray}
for some objects $R_0, R_1 \in \mpi$, and some minimal projective resolution
\begin{eqnarray*} {0\ra P_2' \xra{a_2'} P_1'\xra{a_1'} P_0'\xra {a_0'}A\ra 0.}\end{eqnarray*} \hfill{$\Box$}
\end{lem}

Given an object  $A \in \mathcal{A}$, take a minimal projective resolution
(\ref{minimal resolution}) and consider the corresponding 2-periodic complex
\begin{eqnarray}\label{def of ca}
 \xymatrix{
 &C_A:=P_1\ar[rr]<1mm>^{{\left(\begin{array}{c}a_1\\0\end{array}\right)}}&&P_0\oplus P_2\ar[ll]^{{\left(\begin{array}{cc}0&a_2 \end{array}\right)}}
 }
\end{eqnarray}

Note that any two minimal projective resolutions of $A$ are isomorphic. So the
complex $C_A$ is well-defined up to isomorphism. In addition we have the homological  groups  $H_0(C_A)=A, H_1(C_A)=0$,
and in $K(\mathcal{A})$ we have $\hat{C_A}=\hat{P_0}+\hat{P_2}-\hat{P_1}=\hat{A}.$

\begin{prop}\label{epi}
Given objects $A_1,A_2 \in \mathcal{A}$,
 take minimal projective resolutions
$$\ 0 \xra{} P_2 \xra{a_2} P_1\xra{a_1} P_0\xra {a_0}A_1\xra{} 0$$
$$\ 0 \xra{} Q_2 \xra{b_2} Q_1\xra{b_1} Q_0\xra {b_0}A_2\xra{} 0$$
then there is surjective $k$-linear map  $\phi: \Hom_{\mc_2(\mathcal{A})}(C_{A_1},C_{A_2})\twoheadrightarrow \Hom_{\mathcal{A}}(A_1,A_2)$
and
 $$|\ker \phi|=|\Hom_{\mathcal{A}}(P_2,Q_0)|\times|\Hom_{\mathcal{A}}(P_1,Q_2)|\times|\Hom_{\mathcal{A}}(P_0,Q_1)|/|\Hom_{\mathcal{A}}(P_0,Q_2)|.$$
\end{prop}
{\bf Proof }:
Let
\[\widetilde{f}=\left(\begin{array}{cc}h\ ,&{\left(\begin{array}{cc}a&b \\c&d\end{array}\right) } \end{array}\right)\in  \Hom_{\mc_2(\mathcal{A})}(C_{A_1},C_{A_2}).
\]
We have the  following commutative diagram
$$ \xymatrix{
&P_1\ar[dd]_{h}\ar[rr]<1mm>^{{\left(\begin{array}{c}a_1\\0\end{array}\right)}}&&P_0\oplus P_2\ar[ll]^{{\left(\begin{array}{cc}0&a_2 \end{array}\right)}}\ar[dd]^{{\left(\begin{array}{cc}a&b \\c&d\end{array}\right)}}\\
&&&\\
&Q_1\ar[rr]<0.5mm>^{{\left(\begin{array}{c}b_1\\0\end{array}\right)}}&&Q_0\oplus Q_2\ar[ll]<0.5mm>^{{\left(\begin{array}{cc}0&b_2 \end{array}\right)}} }
$$
i.e., we have
\[
\left\{
\begin{array}{ccc}
a\circ a_1&=&b_1\circ h\\
c\circ a_1&=&0\\
b_2\circ c&=&0\\
b_2\circ d&=&h\circ a_2
\end{array}
\right.
\]
Because $b_2$ is a monomorphism, we have $c=0$. And we have the following commutative diagram
$$ \xymatrix{
0\ar[r]&P_2\ar[r]^{a_2}\ar[d]_d&P_1\ar[r]^{a_1}\ar[d]_h&P_0\ar[r]^{a_0}\ar[d]_a&A_1\ar[r]&0\\0\ar[r]&Q_2\ar[r]^{b_2}&Q_1\ar[r]^{b_1}&Q_0\ar[r]^{b_0}&A_2\ar[r]&0
}$$
So there exists a unique $f\in \Hom_\mathcal{A}(A_1,A_2)$ such that the diagram
$$ \xymatrix{
0\ar[r]&P_2\ar[r]^{a_2}\ar[d]_d&P_1\ar[r]^{a_1}\ar[d]_h&P_0\ar[r]^{a_0}\ar[d]_a&A_1\ar[r]\ar@{.>}[d]^f&0\\0\ar[r]&Q_2\ar[r]^{b_2}&Q_1\ar[r]^{b_1}&Q_0\ar[r]^{b_0}&A_2\ar[r]&0
}$$
commutes. We define $\phi(\widetilde{f})=f.$

On the other hand,such that the diagram  $$ \xymatrix{
0\ar[r]&P_2\ar[r]^{a_2}\ar@{.>}[d]_d&P_1\ar[r]^{a_1}\ar@{.>}[d]_h&P_0\ar[r]^{a_0}\ar@{.>}[d]_a&A_1\ar[r]\ar[d]^f&0\\0\ar[r]&Q_2\ar[r]^{b_2}&Q_1\ar[r]^{b_1}&Q_0\ar[r]^{b_0}&A_2\ar[r]&0
}$$
commutes. For any $b\in\Hom_{\ma}(P_2,Q_0)$, we have
 \[\widetilde{f}=\left(\begin{array}{cc}h\ ,&{\left(\begin{array}{cc}a&b \\0&d\end{array}\right) } \end{array}\right)\in\Hom_{\mc_2(\ma)}(C_{A_1},C_{A_2}),\]
 and $\phi(\widetilde{f})=f$ by the definition.
Therefore,
 there is a surjection   $$\phi: \Hom_{\mc_2(\mathcal{A})}(C_{A_1},C_{A_2})\twoheadrightarrow \Hom_{\mathcal{A}}(A_1,A_2).$$ It is easy to see that $\phi$  is a $k$-linear map.

Let  \[\widetilde{f}=\left(\begin{array}{cc}h\ ,&{\left(\begin{array}{cc}a&b \\c&d\end{array}\right) } \end{array}\right)\in\Hom_{\mc_2(\ma)}(C_{A_1},C_{A_2}),\] such that $\phi(\widetilde{f})=0$. We have $c=0$ and
the following exchange graph $$ \xymatrix{
0\ar[r]&P_2\ar[r]^{a_2}\ar[d]_d&P_1\ar[r]^{a_1}\ar[d]_h&P_0\ar[r]^{a_0}\ar[d]_a&A_1\ar[r]\ar[d]^0&0\\0\ar[r]&Q_2\ar[r]^{b_2}&Q_1\ar[r]^{b_1}&Q_0\ar[r]^{b_0}&A_2\ar[r]&0
}.$$
So there exists only a morphism $a'$ from $P_0$ to $\ker(b_0)$. Conversely, consider the exact sequence
$$\ 0 \ra \ker b_0 \xra{i_0} Q_0\xra{b_0} A_2 \ra 0$$
and if there is a morphism $a'$ from
$P_0$ to $\ker(b_0)$, we can let $a=a'\circ i_0$. So the number of $a\in\Hom_{\ma}(P_0,Q_0)$ with $\phi(\widetilde{f})=0$ is equal to $|\Hom_{\ma}(P_0,\ker(b_0))|.$
Since $P_0$ is a projective module,
using the exact sequence
 $$\ 0 \ra Q_2 \xra{b_2} Q_1\xra{\pi_1} \ker(b_0) \ra 0$$
we have  $\dim_k \Hom_{\mathcal{A}}(P_0, \ker(b_0))=\dim_k \Hom_{\mathcal{A}}(P_0,Q_1)-\dim_k \Hom_{\mathcal{A}}(P_0,Q_2),$ i.e. $$|\Hom_{\mathcal{A}}(P_0, \ker(b_0))|=|\Hom_{\mathcal{A}}(P_0,Q_1)|/|\Hom_{\mathcal{A}}(P_0,Q_2)|.$$

Now fix $a$. Since   $d$ is determined uniquely by $h$,  we just need to determine $h$.
Note that $h$  satisfies the following  exchange graph
\begin{eqnarray}\label{exchange graph1} \xymatrix{
P_1\ar[r]^{a_1}\ar[d]_h&P_0\ar[d]_{a}\\Q_1\ar[r]^{b_1}&Q_0
}\end{eqnarray}
This reduces the  following exchange graph
\begin{eqnarray}\label{exchange graph2} \xymatrix{
P_1\ar[r]^{a_1}\ar[d]_h&P_0\ar[d]_{a'}\\Q_1\ar@{>>}[r]^{\pi_1}&\ker b_0
}\end{eqnarray}
Since $P_1$ is projective,  using the exact sequence
$$\ 0 \ra Q_2 \xra{b_2} Q_1\xra{\pi_1} \ker b_0 \ra 0,$$
we get an exact sequence
$$0\rightarrow \Hom_{\mathcal{A}}(P_1,Q_2)\rightarrow  \Hom_{\mathcal{A}}(P_1,Q_1)\rightarrow \Hom_{\mathcal{A}}(P_1, \ker(b_0)) \rightarrow 0 .$$
Note that $a' \circ a_1\in \Hom(P_1,\ker(b_0))$ is fixed. So the number of $h$ satisfying the  exchange graph (\ref{exchange graph2}) is equal  to $|\Hom_{\mathcal{A}}(P_1,Q_2)|$.

From the  above discussions and noting that  there is no restriction to the morphism $b$ from $P_2$ to $Q_0$, we have  $$|\ker(\phi)|=|\Hom_{\mathcal{A}}(P_2,Q_0)|\times|\Hom_{\mathcal{A}}(P_1,Q_2)|\times|\Hom_{\mathcal{A}}(P_0,Q_1)|/|\Hom_{\mathcal{A}}(P_0,Q_2)|.$$
 \hfill{$\Box$}

 By Lemma \ref{rootcat and homolic cat} and Lemma \ref{ext and hom over 2-}, we have the following isomorphisms:
\begin{lem}\label{ext iso}
 $\Ext_{\mc_2(\mathcal{A})}(C_{A_1},C_{A_2})\cong\Ext_{\mathcal{A}}(A_1,A_2).$
\end{lem}
 {\bf Proof }:
\begin{eqnarray*}\Ext_{\mc_2(\mathcal{A})}(C_{A_1},C_{A_2})&\cong&\Hom_{\mk_2(\mathcal{A})}(C_{A_1},C_{A_2}^*)\\&\cong &\Hom_{R_2({\ma})}({A_1},{A_2}[1])\\&\cong&\bigoplus_{i\in \mathbb{Z}}\Hom_{\md^b(\ma)}(A_1,A_2[2i+1])\\&=& \bigoplus_{i\geq 0}\Hom_{\md^b(\ma)}(A_1,A_2[2i+1]) \\&=&\bigoplus_{i\geq 0}\Ext_{\mathcal{A}}^{2i+1}(A_1,A_2)
\end{eqnarray*}
Since the global dimension of $\ma$ is at most two, we have $$\Ext_{\mc_2(\mathcal{A})}(C_{A_1},C_{A_2})\cong\bigoplus_{i\geq 0}\Ext_{\mathcal{A}}^{2i+1}(A_1,A_2)=\Ext_{\mathcal{A}}(A_1,A_2).$$
 \hfill{$\Box$}

Given $A\in \ma$, take a minimal projective resolution (\ref{minimal resolution}), define
$${E_A}=t^{\langle \hat{P_1}-2\hat{P_2}, \hat{A}\rangle}\cdot K_{-\hat{P_1}}\ast  K_{\hat{P_2}} \ast  K_{-\hat{P_2}}^* \ast [C_A]\in \md\mh(\ma).$$
From the following proposition, one can see that $E_A$ is not depend  on whether the projective resolution is minimal or not.
\begin{prop}
Suppose we  take a different,
not necessarily minimal, projective resolution (\ref{not minimal resolution}), consider the corresponding complex
(\ref{def of ca}) in $\mc_2(\mathcal{P})$, i.e. let $$ \xymatrix{
&C_A':=R_0\oplus P_1\oplus R_1\ar[rrr]<1mm>^{{\left(\begin{array}{ccc}1&0&0\\0&a_1&0\\0&0&0\\0&0&0\end{array}\right)}}&&&R_0\oplus P_0\oplus R_1\oplus P_2\ar[lll]^{{\left(\begin{array}{cccc}0&0&0&0\\0&0&0&a_2\\0&0&1&0 \end{array}\right)}}
 }, $$\noindent and  $${E_A'}=t^{\langle \hat{P_1}+ \hat{R_0}+\hat{R_1}-2(\hat{P_2}+\hat{R_1}), \hat{A}\rangle}\cdot K_{-\hat{P_1}-\hat{R_0}-\hat{R_1}}\ast  K_{\hat{P_2}+\hat{R_1}} \ast  K_{-\hat{P_2}-\hat{R_1}}^* \ast [C_A'].$$
 \noindent Then we have
 $${E_A'}={E_A}.$$
\end{prop}
 {\bf Proof }:
Note that $[C_A']=[ K_{R_1}^* \oplus K_{R_0} \oplus  C_A]$. By Lemma \ref{fact gongshi} and Lemma \ref{fact gongshi2},
we have
\begin{eqnarray*}
[C_A']&=&[ K_{R_1}^* \oplus K_{R_0} \oplus  C_A]
=t^{\langle \hat{R_1},\hat{A}\rangle-\langle \hat{R_0}, \hat{A}\rangle}\cdot K_{\hat{R_1}}^*\ast K_{\hat{R_0}}\ast [C_A].
\end{eqnarray*}
Hence,
\begin{eqnarray*}
  {E_A'}&=& t^{\langle \hat{P_1}+ \hat{R_0}+\hat{R_1}-2(\hat{P_2}+\hat{R_1}), \hat{A}\rangle}\cdot K_{-\hat{P_1}-\hat{R_0}-\hat{R_1}}\ast  K_{\hat{P_2}+\hat{R_1}} \ast  K_{-\hat{P_2}-\hat{R_1}}^* \ast [C_A']\\
 &=& t^{\langle \hat{P_1}+ \hat{R_0}+\hat{R_1}-2(\hat{P_2}+\hat{R_1}), \hat{A}\rangle}\cdot K_{-\hat{P_1}-\hat{R_0}-\hat{R_1}}\ast  K_{\hat{P_2}+\hat{R_1}} \ast  K_{-\hat{P_2}-\hat{R_1}}^*\notag\\& & \ast t^{\langle \hat{R_1},\hat{A}\rangle-\langle \hat{R_0}, \hat{A}\rangle}\cdot K_{\hat{R_1}}^*\ast K_{\hat{R_0}}\ast [C_A]\\
  &=& t^{\langle \hat{P_1}+ \hat{R_0}+\hat{R_1}-2(\hat{P_2}+\hat{R_1}), \hat{A}\rangle+\langle \hat{R_1},\hat{A}\rangle-\langle \hat{R_0}, \hat{A}\rangle}\cdot K_{-\hat{P_1}-\hat{R_0}-\hat{R_1}}\ast  K_{\hat{P_2}+\hat{R_1}} \ast  K_{-\hat{P_2}-\hat{R_1}}^*\notag\\& & \ast  K_{\hat{R_1}}^*\ast K_{\hat{R_0}}\ast [C_A]\\
  &=& t^{\langle \hat{P_1}-2\hat{P_2}, \hat{A}\rangle}\cdot K_{-\hat{P_1}}\ast  K_{\hat{P_2}} \ast  K_{-\hat{P_2}}^* \ast [C_A]\\
  &=& E_A
\end{eqnarray*}
 \hfill{$\Box$}
\begin{prop}\label{linear indepent}
For $B_i\in\ma,\  1\leq i \leq s,$ assume that $\{[B_i]|1\leq i \leq s\}$ is linearly independent in $\mathcal{H}_{tw}(\mathcal{A})$. Then for any $P_i,Q_i,R_i,S_i(1\leq i\leq s)\in \mpi$, $\{K_{\hat{P_i}}\ast K_{\hat{Q_i}}^*\ast K_{-\hat{R_i}}\ast K_{-\hat{S_i}}^*\ast[C_{B_i}]|1\leq i \leq s\}$  is linearly independent in $\mathcal{DH}(\mathcal{A}).$
\end{prop}
{\bf Proof }:
Assume there exist $l_i\in \mathbb{C},\  1\leq i \leq s,$ such that in $\mathcal{DH}(\mathcal{A})$
$$\sum_{i=1}^s l_i K_{\hat{P_i}}\ast K_{\hat{Q_i}}^*\ast K_{-\hat{R_i}}\ast K_{-\hat{S_i}}^*\ast[C_{B_i}]=0.$$
 Then multiplying  suitable $K_{\hat{P}}, K_{\hat{Q}}$ for some $P,Q\in \mpi$ we can get
 $$\sum_{i=1}^s l_i K_{\hat{X_i}}\ast K_{\hat{Y_i}}^*\ast[C_{B_i}]=0$$
in $\mh_{tw}(\mc_2(\mpi))$ for some $X_i,Y_i\in \mpi,\  1\leq i \leq s$.
 By Lemma \ref{fact gongshi} and Lemma \ref{fact gongshi2}, we have $$\sum_{i=1}^s l_i\cdot t^{\langle \hat{X_i}-\hat{Y_i}, \hat{B_i} \rangle } \cdot [K_{{X_i}}\oplus K_{{Y_i}}^*\oplus C_{B_i}]=0$$ in $\mh_{tw}(\mc_2(\mpi))$.
If there is $u\ (1\leq u\leq s)$ such that  $l_u\neq 0$, there must be $v\ (1\leq v\leq s)$ such that $v\neq u, l_v\neq 0$ and
$$ K_{{X_{u}}}\oplus K^*_{{Y_{u}}}\oplus C_{B_u}\cong K_{{X_{v}}}\oplus K^*_{{Y_{v}}}\oplus C_{B_v}.$$
Then
$$H_0( K_{{X_{u}}}\oplus K^*_{{Y_{u}}}\oplus C_{B_u})\cong H_0(K_{{X_{v}}}\oplus K^*_{{Y_{v}}}\oplus C_{B_v}).
$$
 Since for any $P\in \mpi$, $H_0(K_P)=H_1(K_P)=H_0(K^*_P)=H_1(K^*_P)=0,$ and $H_0(C_{B_i})=B_i\ \ (1\leq i\leq s)$,
 we have  $B_u\cong B_v,$   a contradiction  with $\{[B_i]|1\leq i \leq s\}$ linearly independent in $\mathcal{H}_{tw}(\mathcal{A})$. So for any $u\ (1\leq u\leq s)$, we have $l_u=0.$ Hence, $\{K_{\hat{P_i}}\ast K_{\hat{Q_i}}^*\ast K_{-\hat{R_i}}\ast K_{-\hat{S_i}}^*\ast[C_{B_i}]|1\leq i \leq s\}$  is linearly independent in $\mathcal{DH}(\mathcal{A}).$
\hfill{$\Box$}

From the above proposition, it is easy to get the following result:
\begin{cor}\label{Ebi linear indepent}
For $B_i\in\ma,\  1\leq i \leq s,$ assume that $\{[B_i]|1\leq i \leq s\}$ is linearly independent in $\mathcal{H}_{tw}(\mathcal{A})$. Then $\{[C_{B_i}]|1\leq i \leq s\}$  is linearly independent in $\mathcal{DH}(\mathcal{A}),$ and $\{[E_{B_i}]|1\leq i \leq s\}$  are also linearly independent in $\mathcal{DH}(\mathcal{A}).$
\hfill{$\Box$}
\end{cor}

\begin{thm}\label{main result}
There is an injective algebra homomorphism
$$I_+: \mathcal{H}_{tw}(\mathcal{A}) \hookrightarrow \mathcal{DH}(\mathcal{A}),\ \ [A] \mapsto {E_A}.$$

\end{thm}
{\bf Proof }:
Given objects $A_1,A_2 \in \mathcal{A}$,
in $\mathcal{H}_{tw}(\mathcal{A})$   we have
\begin{eqnarray*}
[{A_1}]\ast[{A_2}]&=&t^{\langle \hat{A_1},\hat{A_2}\rangle}\sum_{[A_3]}\frac{|\Ext_{\mathcal{A}}({A_1},{A_2})_{A_3}|}{|\Hom_{\mathcal{A}}({A_1},{A_2})|}\cdot [{A_3}].
\end{eqnarray*}
Take minimal projective resolutions
$$\ 0 \ra P_2 \xra{a_2} P_1\xra{a_1} P_0\xra {a_0}A_1\ra 0$$
$$\ 0 \ra Q_2 \xra{b_2} Q_1\xra{b_1} Q_0\xra {b_0}A_2\ra 0.$$
Applying the Horseshoe Lemma gives the projective resolution of ${A_3}$,
$$\ 0 \ra P_2\oplus Q_2 \xra{} P_1\oplus Q_1\xra{} P_0\oplus Q_0\xra {}{A_3}\ra 0.$$

 By Lemma \ref{ext iso}, $\Ext_{\mc_2(\mathcal{A})}(C_{A_1},C_{A_2})=\Hom_{\mk_2(\mathcal{A})}(C_{A_1},C_{A_2}^*)=\Ext_{\mathcal{A}}(A_1,A_2)$. Then it is easy to see that any extension of $C_{A_1}$ by $C_{A_2}$ is the complex $C_{A_3}$ defined by the corresponding extension $A_3$ of $A_1$ by $A_2$.
Using Lemma \ref{fact gongshi} and Lemma \ref{fact gongshi2},
we have
\begin{eqnarray}
& &I_+({A_1})\ast I_+({A_2})\notag\\
&=&t^{\langle \hat{P_1}-2\hat{P_2},\hat{A_1}\rangle+\langle \hat{Q_1}-2\hat{Q_2},\hat{A_1}\rangle}\cdot K_{-\hat{P_1}}\ast K_{\hat{P_2}}\ast K_{-\hat{P_2}}^*\ast[C_{A_1}]\ast K_{-\hat{Q_1}}\ast K_{\hat{Q_2}}\ast K_{-\hat{Q_2}}^*\ast[C_{A_2}]\notag\\
&=&t^{\langle \hat{P_1}-2\hat{P_2},\hat{A_1}\rangle+\langle \hat{Q_1}-2\hat{Q_2},\hat{A_2}\rangle+(\hat{Q_1},\hat{A_1})-2(\hat{Q_2},\hat{A_1})}\cdot K_{-\hat{P_1}}\notag\\& &\ast K_{\hat{P_2}}\ast K_{-\hat{P_2}}^*\ast K_{-\hat{Q_1}}\ast K_{\hat{Q_2}}\ast K_{-\hat{Q_2}}^*\ast[C_{A_1}]\ast[C_{A_2}]\notag\\
&=&t^{\langle \hat{P_1}-2\hat{P_2},\hat{A_1}\rangle+\langle \hat{Q_1}-2\hat{Q_2},\hat{A_2}\rangle+(\hat{Q_1},\hat{A})-2(\hat{Q_2},\hat{A})+\langle \hat{P_1},\hat{Q_1}\rangle+\langle \hat{P_0}+\hat{P_2},\hat{Q_0}+\hat{Q_2}\rangle-\langle \hat{P_1}+\hat{Q_1}-2\hat{P_2}-2\hat{Q_2},\hat{{A_3}}\rangle}\notag\\
& & \cdot\sum_{[{A_3}]}\frac{|\Ext_{\mathcal{A}}({A_1},{A_2})_{A_3}|}{|\Hom_{\mc_2(\mathcal{A})}(C_{A_1},C_{A_2})|}\ast E_{A_3}.\label{result2}
\end{eqnarray}
From Proposition \ref{epi}, we know
\begin{eqnarray*}
 & &|\Hom_{\mc_2(\mathcal{A})}(C_{A_1},C_{A_2})|\\&=&|\Hom_{\mathcal{A}}(A_1,A_2)|\cdot|\Hom_{\mathcal{A}}(P_2,Q_0)|\cdot|\Hom_{\mathcal{A}}(P_1,Q_2)|
 \cdot|\Hom_{\mathcal{A}}(P_0,Q_1)|/|\Hom_{\mathcal{A}}(P_0,Q_2)|\end{eqnarray*}

 Because
\begin{eqnarray}
& &I_+([{A_1}]\ast[{A_2}])=t^{\langle \hat{A_1},\hat{A_2}\rangle}\cdot\sum_{{[A_3]}}\frac{|\Ext_{\mathcal{A}}({A_1},{A_2})_{A_3}|}{|\Hom_{\mathcal{A}}({A_1},{A_2})|}\ast E_{A_3},
\label{result1}
\end{eqnarray}
comparing (\ref{result2}) and (\ref{result1}), since $|\Hom_{\mathcal{A}}(P,Q)|=q^{\langle \hat{P},\hat{Q}\rangle}=t^{2\langle \hat{P},\hat{Q}\rangle}$ for any projective modules $P,Q$, we have $I_+([{A_1}]\ast[{A_2}])=I_+([{A_1}])\ast I_+([{A_2}])$ if and only  if
\begin{eqnarray}
& &t^{\langle \hat{P_1}-2\hat{P_2},\hat{A_1}\rangle+\langle \hat{Q_1}-2\hat{Q_2},\hat{A_2}\rangle+(\hat{Q_1},\hat{A_1})-2(\hat{Q_2},\hat{A_1})+\langle \hat{P_1}, \hat{Q_1}\rangle+\langle \hat{P_0}+\hat{P_2},\hat{Q_0}+\hat{Q_2}\rangle-\langle \hat{P_1}+\hat{Q_1}-2\hat{P_2}-2\hat{Q_2},\hat{{A_3}}\rangle}\notag\\
&=&t^{\langle \hat{A_1},\hat{A_2}\rangle+2\langle \hat{P_2},\hat{Q_0}\rangle+2\langle \hat{P_1},\hat{Q_2}\rangle+2\langle \hat{P_0},\hat{Q_1}\rangle-2\langle \hat{P_0},\hat{Q_2}\rangle}\label{zhishuxiangdeng}
\end{eqnarray}
Using $\hat{A_3}=\hat{A_1}+\hat{A_2},\hat{A_1}=\hat{P_0}+\hat{P_2}-\hat{P_1}$ and $\hat{A_2}=\hat{Q_0}+\hat{Q_2}-\hat{Q_1}$, it is easy to verify that the equality (\ref{zhishuxiangdeng}) holds .

By Corollary \ref{Ebi linear indepent}, if $\{[B_i]|1\leq i\leq s\}$ is linearly independent in $\mathcal{H}_{tw}(\mathcal{A})$ for any $B_i(1\leq i\leq s)\in \ma$,
then  $\{[E_{B_i}]|1\leq i\leq s\}$ is also linearly independent in $\mathcal{DH}(\mathcal{A})$. Therefore, $I_+$ is injective.
\hfill{$\Box$}

\

Let $\pi$ be the natural homomorphism from $\mathcal{DH}(\mathcal{A})$ to $\mathcal{DH}_{red}(\mathcal{A}),$  we have
\begin{thm}\label{reduce embedding}
There is an embedding of algebras $\widetilde{I_+}=\pi \circ I_+ : \mh_{tw}(\mathcal{A}) \hookrightarrow \mathcal{DH}_{red}(\mathcal{A})$.
\end{thm}
{\bf Proof }: Clearly $\widetilde{I_+}$ is an algebra homomorphism. So we only need to prove that $\widetilde{I_+}$ is injective.
Let $J$ be the ideal in $\mathcal{DH}(\mathcal{A})$ generated by $K_{-\hat{P}}-K^*_{\hat{P}}$ for all $P\in \mpi$. Then by above theorem and  $\mathcal{DH}_{red}(\mathcal{A})=\mathcal{DH}(\mathcal{A})/J$, we just need to prove $\mathrm{im} I_+\cap J=0$.

It is easy to see that all $K_{\hat{P}}\ast K^*_{\hat{Q}}\ast K_{-\hat{R}}\ast K^*_{-\hat{S}}\ast [M_{\bullet}]$ span $\mathcal{DH}(\mathcal{A})$, where
$P, Q, R, S \in \mpi$ and $ M_{\bullet}\in \mc_2(\mpi)$ has no non-zero direct summand of acyclic complexes.
Therefore any element in $\mathrm{im} I_+\cap J$ has the form
\begin{eqnarray*} \sum_{i=1}^sl_iE_{B_i}=\sum_{j=1}^nh_j( K_{-\hat{X_j}}-K^*_{\hat{X_j}})\ast K_{\hat{T_j}}\ast K^*_{\hat{W_j}}\ast K_{-\hat{Y_j}}\ast K^*_{-\hat{Z_j}}\ast[(M_{j})_{\bullet}]
\label{reduce emddng eqn}\end{eqnarray*}
for some $l_i,h_j\in \mathbb{C},\ 1\leq i\leq s,\ 1\leq j\leq n,$ and some  $B_i\in \ma, \ 1\leq i\leq s$, where  $\{[B_i]|1\leq i\leq s\}$ is linearly independent in $\mathcal{H}_{tw}(\mathcal{A})$, and some $T_j, W_j, X_j,Y_j,Z_j\in \mpi$ and $(M_{j})_{\bullet}\in \mc_2(\mpi)$ for $ 1\leq j\leq n$ such that each $(M_{j})_{\bullet}$ has no non-zero direct summand of acyclic complexes.
Multiplying $K_P, K^*_Q$ on the both sides for some suitable $P,Q\in \mpi$, we can get
$$\sum_{i=1}^s l_i t_i K_{\hat{P_{i}}}\ast K^*_{\hat{Q_{i}}}\ast [C_{B_i}]=\sum_{j=1}^nh_j(1- K_{\hat{X_j}}\ast K^*_{\hat{X_j}})\ast K_{\hat{R_{j}}}\ast K^*_{\hat{S_{j}}}\ast [(M_{j})_{\bullet}]$$
in $\mathcal{H}_{tw}({\mc_2(\mpi)})$
for some $P_i,Q_i,R_j,S_j\in \mpi,\ 1\leq i\leq s,1\leq j\leq n$, where all $t_i\in \mathbb{C}$ are not zero.

Suppose some $l_i\neq 0$ and take $(M_{i1})_{\bullet}, (M_{i2})_{\bullet},\cdots, $ $(M_{im_i})_{\bullet}$ $\in \{(M_1)_{\bullet},(M_2)_{\bullet},\cdots, (M_n)_{\bullet}\}$ to be all the $\mathbb{Z}_2$-graded complexes such that each of them is isomorphic to  $C_{B_i}$. Note that
$$\{\ K_{\hat{P}}\ast K^*_{\hat{Q}}\ast [M_{\bullet}]\ | \ P, Q \in \mpi, \ M_{\bullet}\in \mc_2(\mpi)\  \mbox{has no non-zero direct summand of acyclic complexes}\}$$
is also a basis in $\mathcal{H}_{tw}({\mc_2(\mpi)}).$ We have
 \begin{eqnarray}\label{reduce emddng eqn2}l_i t_i K_{\hat{P_{i}}}\ast K^*_{\hat{Q_{i}}}\ast [C_{B_i}]=\sum_{j=1}^{m_i} h_{ij}(1- K_{\hat{X_{ij}}}\ast K^*_{\hat{X_{ij}}})\ast K_{\hat{R_{{ij}}}}\ast K^*_{\hat{S_{{ij}}}}\ast [(M_{ij})_{\bullet}]
 \end{eqnarray}
in $\mathcal{H}_{tw}({\mc_2(\mpi)})$.
Under the above basis,  since the sum of the coefficients of the equation (\ref{reduce emddng eqn2}) on the right is zero, we have $l_it_i=0$ as the  coefficient on the left and so $l_i=0$,  which is  a contradiction  with $l_i \neq 0$.
Hence for any $1\leq i\leq s$, we have $l_i=0$. Therefore, $\mathrm{im} I_+\cap J=0.$
\hfill{$\Box$}

\begin{rks}
If we set $F_A=E_A^*$, then there is also an injective algebra homomorphism $$I_-: \mathcal{H}_{tw}(\mathcal{A}) \hookrightarrow \mathcal{DH}(\mathcal{A}),\ \ [A] \mapsto {F_A}$$ and  this  embedding  can induce an embedding  form  $ H_{tw}(\mathcal{A})$ to $\mathcal{DH}_{red}(\mathcal{A})$.\hfill{$\Box$}
\end{rks}

In the following we will get some interesting consequences.

Let us recall the definitions of tilted algebras and canonical algebras.

\begin{df}$\mathrm{(\cite{HR})}$ Let $H$ be a finite-dimensional hereditary algebra. An $H$-module $T$ is a tilting module if it satisfies:

 $(1)$ proj.dim$_H$T$\leq 1$;

$(2)$ $\Ext_H^1(T,T)=0$;

$(3)$ there exists a short exact sequence $0\ra H\ra T_0 \ra T_1\ra 0$, where $T_0,T_1$ are direct summands of a finite direct sum of copies of $T$ respectively.

The endomorphism algebra $B=\End_HT$ is called a tilted algebra.
\end{df}

\begin{df}$\mathrm{(\cite{Ri1})}$
Let $\emph{\textbf{p}}=(p_0,p_1,\cdots,p_n)$ be a $(n+1)$-tuple of integers $p_i\geq 1$, $\bm{\lambda}=(\lambda_2,\cdots,\lambda_n)$ be  pairwise distinct elements of $k$. A canonical algebra $B$ of type $(\emph{\textbf{p}},\bm{\lambda})$ is the following quiver

$$ 
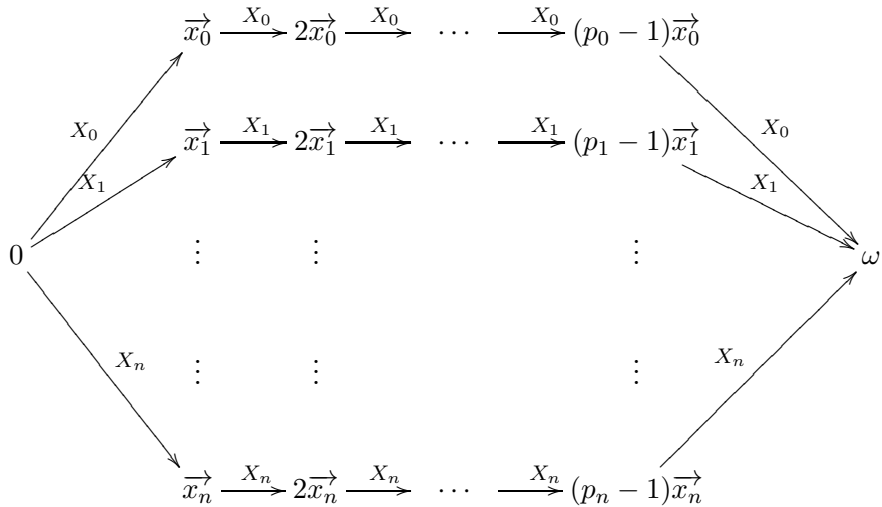
\xymatrix{
&&\overrightarrow{x_0}\ar[r]^{X_0} &2\overrightarrow{x_0}\ar[r]^{X_0} &\ \ \cdots\ \ \ar[r]^{X_0}&(p_0-1)\overrightarrow{x_0}\ar[ddrr]^{X_0} &&\\
&&\overrightarrow{x_1}\ar[r]^{X_1} &2\overrightarrow{x_1}\ar[r]^{X_1} &\ \ \cdots\ \ \ar[r]^{X_1}&(p_1-1)\overrightarrow{x_1}\ar[drr]^{X_1} &&
 \\  0\ar[urr]^{X_1}\ar[uurr]^{X_0}\ar[ddrr]^{X_n} & &\vdots&\vdots&&\vdots& &   \omega       \\&&\vdots&\vdots&&\vdots&&
 \\&&\overrightarrow{x_n}\ar[r]^{X_n} &2\overrightarrow{x_n}\ar[r]^{X_n} &\ \ \cdots\ \ \ar[r]^{X_n}&(p_n-1)\overrightarrow{x_n}\ar[uurr]^{X_n} &&
}$$
with relations given by
$$X_i^{p_i}= X_0^{p_0}-\lambda_iX_1^{p_1} \ \ \  \mathrm{for}\ \ \  i = 2,3,\ldots,n, $$
 By \cite{GL}, $\md^b(B)$ is equivalent to $\md^b(coh(\mathbb{X}))$, where $coh(\mathbb{X})$ is the category of coherent sheaves on a weighted projective line $\mathbb{X}$ of type $(\emph{\textbf{p}},\bm{\lambda)}$.
\end{df}

 By \cite{HR} tilted algebras have global dimension at most two and it is obvious that canonical algebras have global dimension $\leq 2.$
As an  application, we have
\begin{cor}
Assume that $\mathcal{A}=\mathrm{mod} B$ is  the category of finitely generated right $B$-modules
 where $B$ is a tilted algebra or a canonical algebra. Then there is an embedding of algebras from $\mh_{tw}(\mathcal{A})$ to $\mathcal{DH}(\mathcal{A})$ and $\mathcal{DH}_{red}(\mathcal{A})$.\hfill{$\Box$}
\end{cor}
\begin{rks} It should be also interesting to apply our main theorem to some other  algebras  which are neither tilted algebras nor tubular algebras but with global dimension $\leq 2$.
For example, the algebra is given respectively by the following quivers with relations
$$  \xymatrix{
\circ\ar@/^/[r]^{\alpha}&\circ\ar@/^/[l]^{\beta}&\alpha\beta=0,
}
$$
and
$$  \xymatrix{&\circ\ar[dr]^{\beta}&&\\
\circ\ar[ur]^{\alpha}\ar[rr]^{\gamma}&&\circ&\beta\alpha=0
}
$$
\end{rks}

\end{document}